\documentclass[12pt]{article}

\title{Generalized Trigonometric Functions over Associative Algebras}

\author{Nathan BeDell \\ nbedell@tulane.edu}

\usepackage[a4paper, total={6in, 8in}]{geometry}
\usepackage{graphicx}
\usepackage{amsmath}
\usepackage{amssymb}
\usepackage{mathrsfs}
\usepackage{amsfonts}
\usepackage[utf8]{inputenc}
\usepackage[english]{babel}
\usepackage{amsthm}
\usepackage{amscd}
\usepackage{tikz}
\usepackage{bbm}

\newcommand{\Acal}{\mathcal{A}}

\newcommand{\Ccal}{\mathcal{C}}
\newcommand{\Hcal}{\mathcal{H}}

\newcommand{\Ld}{\mathrm{Ld}}

\newcommand{\comp}{ \,{\scriptstyle \stackrel{\circ}{}}\, }

\theoremstyle{definition}
\newtheorem{definition}{Definition}[section]

\newtheorem{theorem}[definition]{Theorem}
\newtheorem{example}[definition]{Example}
\newtheorem{lemma}[definition]{Lemma}
\newtheorem{corollary}[definition]{Corollary}
\newtheorem{proposition}[definition]{Proposition}

\theoremstyle{remark}
\newtheorem*{open}{Open Problem}

\theoremstyle{remark}

\theoremstyle{remark}

\theoremstyle{remark}

\begin{document}

\maketitle

\begin{abstract}
Extending the work of Freese \cite{Freese15}, we further develop the theory of generalized trigonometric functions. In particular, we study to what extent the notion of polar form for the complex numbers may be generalized to arbitrary associative algebras, and how the general trigonometric functions may be used to give particularly elegant formulas for the logarithm over an algebra. Finally, we close with an array of open questions relating to this line of inquiry.
\end{abstract}

\section{Introduction to $\Acal$-Calculus}
\label{sec:intro_a_cal}

Mathematicians have long considered generalizations of the usual complex analysis over different algebraic context, for instance, the hyperbolic numbers \cite{LorentzNumbers} are well studied, and Price has studied what in our context we call ``type-C algebras'' in some generality \cite{Price} \footnote{For a more complete review of the relevant literature, and how our approach relates, see \cite{Cook12}.}. In this paper, we work in the context of unital associative commutative finite dimensional algebras, which henceforth we shall just call ``algebras". Although conceivably our work could be generalized, working in this context gives us a particularly nice theory of differentiation, power series, as well as exponential and logarithmic functions over an algebra, as was originally studied in \cite{Cook12}, \cite{Freese15}, and \cite{algebrapaper}. As we will see later, this theory which has been built up will be instrumental in our exposition of the generalized trigonometric functions over an algebra. 

A fuller exposition of our ``$\Acal$-calculus'' program, developing the theory of basic calculus, analysis, and beyond in the context of real associative algebras can be found in the aforementioned papers, but for the sake of making this paper self-contained, we will give a brief introduction to the topic here.

\begin{definition}
	An algebra $\mathcal{A}$ is a finite dimensional real vector space together with a bilinear multiplication operation $\star : \mathcal{A} \times \mathcal{A} \rightarrow \mathcal{A}$ satisfying the following properties:
	\begin{enumerate}
		\item $ v \star (w \star z) = (v \star w) \star z$ for all $v,w,z \in \mathcal{A}$
		\item There exists an element $\mathbbm{1} \in \mathcal{A}$ such that $\mathbbm{1} \star z = z \star \mathbbm{1} = z$ for all $z \in \mathcal{A}$.
	\end{enumerate}
\end{definition}

The basic idea in $\Acal$-calculus is to exploit the vector space structure inherent in an algebra to define the operation of differentiation in terms of the usual Frechet differential from advanced calculus. Specifically, this is done by way of the matrix representation of algebra elements:

\begin{definition}
	Given an algebra $\mathcal{A}$ the set of all linear transformations $T: \Acal \rightarrow \Acal$ for which $T(x \star y) = T(x) \star y$ forms the \textit{regular representation} which we denote $\mathcal{R}_{\Acal}$. Clearly $T(x) = T(1 \star x) = T(1) \star x$ hence the regular representation is formed by left-multiplications of $\Acal$. 
	
	Denote the left-multiplication by $\alpha \in \Acal$ by $L_\alpha(x) = \alpha \star x$ for each $x \in \Acal$. Given basis $\beta = \{ v_1, \dots, v_n \}$, 
	we define $M_\beta(\alpha) = [L_\alpha]_\beta$ and denote the collection of all such matrices by $M_\beta(\mathcal{A})$. 
\end{definition}

It is important to understand how to calculate with this definition, so we illustrate an example using the familiar case of $\mathbb{C}$:

\begin{example}\label{ex:2}
	Let $\beta = \{ 1, i \}$ and let $x = x + i y$ be an arbitrary element of $\mathbb{C}$. Applying the same method as we did in the preceding example, we find $(x + i y) \star 1 = x + i y$ and $(x + i y) \star i = -y + i x$. Hence: $$ M_\beta(z) = \begin{bmatrix} x & -y \\ y & x \end{bmatrix} $$
\end{example}

\noindent With this, we are able to define differentiation in the manner already mentioned:

\begin{definition}
	Given a function $f : \mathcal{A} \rightarrow \mathcal{A}$, we say that $f$ is $\mathcal{A}$-differentiable at $z_0$ if the Frechet differential $df_{z_0}$ exists, and $df_{z_o} \in \mathcal{R}_{\Acal}$. 
\end{definition}

This satisfies all of the usual properties one might be familiar with in the usual case of calculus over $\mathbb{R}$ or $\mathbb{C}$, but also applies to more general examples, some of which we will introduce in the next section. 

\begin{proposition}
	\label{prop:derivative_properties}  For $\mathcal{A}$-differentiable functions $f,g : \mathcal{A} \rightarrow \mathcal{A}$,
	\begin{enumerate}
		\item $(f(g(z)))' = f'(g(z)) \star g'(z)$.
		\item $\frac{d}{dz} z^n = n z^n$ for powers $n \in \mathbb{N}$.
		\item $(f(z)\star g(z))' = f'(z) \star g(z) + f(z) \star g'(z)$ (here we assume $\Acal$ is commutative)
		\item if $c \in \mathcal{A}$ a constant then $(c f(z) + g(z))' = c f'(z) + g'(z)$.
	\end{enumerate}
\end{proposition}

The matrix representation of algebra elements will also be important to us later on because of the following observation:

\begin{theorem}
	\label{thm:regular_rep_cols}
	Let $\mathcal{A}$ be an algebra with basis $ \beta = \{ 1,v_2, \dots, v_n \}$.  Then if $u_1, \dots, u_n$ are the columns of $M_\beta(z)$, i.e. $M_\beta(z) = [u_1|u_2|\dots|u_n]$ then $u_1=z$, $u_2 = z \star v_2, \dots , u_n = z \star v_n$. Hence, $M_\beta(z) = [z| z \star v_2 | \dots | z \star v_n]$.
	
	\begin{proof}
		See \cite{Cook12}.
	\end{proof}
\end{theorem}

\subsection{Algebra Presentations}

An important concept we will use to develop the theory of generalized trigonometric functions over an algebra in this paper is that of an algebra presentation. For the purposes of making this a self-contained paper, we present the basic ideas here, but a more in depth account may be found in \cite{logarithms}. By a presentation of an algebra we simply mean an algebra of the form $\mathbb{R}[x_1,x_2,\dots,x_n]/I$ for some ideal I. For example, $\mathbb{R}[i]/\lbrack i^2 + 1 \rbrack$ is a presentation of the complex numbers $\mathbb{C}$.

\begin{definition}[Standard presentations of typical algebras] $ $
	\label{def:common_algebras}
	\begin{enumerate}
		\item The $n$-hyperbolic numbers: $\mathcal{H}_n := \mathbb{R}[j]/\langle j^n - 1 \rangle$
		\item The $n$-complicated numbers: $\mathcal{C}_n := \mathbb{R}[i]/\langle i^n + 1 \rangle$
		\item The $n$-nil numbers: $ \mathbf{\Gamma}_n := \mathbb{R}[\epsilon]/\langle \epsilon^n \rangle$
		\item The total $n$-nil numbers: $\mathbf{\Xi}_n := \mathbb{R}[\epsilon_1,\dots,\epsilon_n]/\langle \epsilon_i\epsilon_j | i,j \in \{ 1, 2, \dots , n \} \rangle$
	\end{enumerate}
\end{definition}

However, some presentations are more useful than others, as some presentations can be degenerate. For instance, $\mathbb{R}[i,j]/\lbrack i^2 - 1, j \rbrack$ is another presentation of $\mathbb{C}$.

\begin{definition}
	Given a presentation $\mathcal{P} = \mathbb{R}[x_1,\dots,x_k]/I$ of an algebra $\mathcal{A}$, we say that the presentation $\mathcal{P}$ is \textit{degenerate} if the set $\{ 1 + I, x_1 + I, \dots , x_k + I \}$ is linearly dependent as a vector space over $\mathbb{R}$. 
\end{definition}

\subsection{Exponential and Logarithmic Functions}

With the theory of differentiation we outlined in section \ref{sec:intro_a_cal}, Freese \cite{Freese15} developed a theory of power series which lets us define the exponential function over an arbitrary algebra.

\begin{definition}
	Given an algebra $\mathcal{A}$, we define the exponential function, denoted either $e^z$ or $\exp(z)$ on an algebra by setting $$ e^z = \sum_{k=0}^{\infty} \frac{z^k}{k!} $$
\end{definition}

In BeDell \cite{logarithms} the basic theory of logarithms over an algebra were studied, in which we discovered that the structure of the logarithms for a certain class of algebras has a particularly nice form, but we first need some preliminary definitions:

\begin{definition}
	In an algebra $\mathcal{A}$, by a \textit{logarithm} we simply mean an inverse to the exponential function $\exp(z) : \mathcal{A} \to \mathcal{A}$. We denote $\mathrm{Im}(\exp) = \Ld(\mathcal{A})$, which we call the \textit{logarithmic domain} for an algebra. Thus, a logarithm on $\mathcal{A}$ is simply a function $\log(z) : X \subseteq \Ld(\mathcal{A}) \to \mathcal{A}$.
\end{definition}

\begin{definition}
	Let $\mathcal{A}$ be an algebra. We say a basis $\beta = \{ v_1, \dots , v_n \}$ of $\mathcal{A}$ is multiplicative if for all $v_i,v_j \in \beta$ we have $v_i \star v_j = c v_k$ for some $c \in \mathbb{R}$ and $v_k \in \beta$. If an algebra $\mathcal{A}$ admits a multiplicative basis, then we say $\mathcal{A}$ is a multiplicative algebra.
\end{definition}

\begin{definition}
	\label{def:typeRC}
	Let $\mathcal{A}$ be a commutative algebra. If there exists another algebra $\mathcal{B}$ such that $\mathcal{A} \cong \mathbb{C} \otimes \mathcal{B}$, then we say $\mathcal{A}$ is a type-C algebra. 
	
	If $\mathcal{A}$ is not isomorphic to the direct product of an algebra $\mathcal{B}$ and some type-C algebra $\mathcal{C}$, then we say that $\mathcal{A}$ is a type-R algebra.
\end{definition}

\noindent From which we obtain the following results\footnote{The derivation of these results involves some complex algebraic machinery, and hence, for the sake of brevity, we omit an exposition here. The reader can read \cite{logarithms} for a more complete treatment.}.

\begin{theorem}
	Every algebra $\mathcal{A}$ has a function $\exp^{-1} 
	: V \subseteq \Ld(\Acal) \rightarrow \mathcal{A}$ which is inverse to the exponential function $\exp : U \subseteq \mathcal{A} \rightarrow \mathcal{A}$.
\end{theorem}

\begin{proposition}
	Let $\Acal$ be a multiplicative algebra. If $\mathcal{A}$ is a type-R algebra, then there exists a unique inverse function to $\exp$ on $\Ld(\Acal)$, denoted $\mathrm{Log}_\mathcal{A}(z)$. Otherwise, there exists infinitely many logarithms determined by the branches of the logarithms defined on the complex portion of the algebra.
\end{proposition}

\begin{proposition}
	\label{prop:type-C-LD}
	For a type-C algebra $\mathcal{A}$, $\mathrm{Ld}(\mathcal{A}) = \mathcal{A}^\times$
\end{proposition}

Thus, to summarize, every multiplicative algebra may be decomposed into the direct product of a type-R and a type-C piece, on which the logarithm behaves much as one would expect on $\mathbb{R}$ and $\mathbb{C}$ respectively. In particular, we obtained the following corollary in \cite{logarithms} for semisimple algebras:

\begin{corollary}
	\label{cor:log_branches_semisimple}
	Given a choice of branch cut for each of the complex components, the logarithm for a semisimple algebra $\mathcal{A}$ may be induced from the product of functions $(\log,\dots ,\log, \mathrm{Log}_1, \dots, \mathrm{Log}_m) : (\mathbb{R}^+)^n \times (\mathbb{C}^\times)^m \rightarrow \mathbb{R}^n \times \mathbb{C}^m$ under the isomorphism $\phi : \mathcal{A} \rightarrow \mathbb{R}^n \times \mathbb{C}^m$ given by Wedderburn's theorem, where $\mathrm{Log}_i : \mathbb{C}^\times \rightarrow \mathbb{C}$ denote the chosen branch cut of the complex logarithm.\footnote{Note that usually a branch cut of the complex logarithm is taken with a domain of $\mathbb{C}^\alpha$ -- that is, $\mathbb{C}$ with a ray starting from 0 removed. This is done to ensure that the logarithm is continuous, but sacrifices having an inverse on all of $\mathrm{Ld}(\mathbb{C}) = \mathbb{C}^\times$, however this corollary can be modified to better suit either approach to branch cuts of the complex logarithm.}
\end{corollary}

However, we did not attempt to describe in this paper how one might go about finding formulas for the logarithms over an algebra analogous to those from complex analysis. For this, we turn to Freese's work on generalized trigonometric functions:

\section{Generalized Trigonometric Functions}

The formula we gave for the logarithm for a semisimple algebra in corollary \ref{cor:log_branches_semisimple} is certainly a straightforward construction, yet if we are working in a different presentation of the algebra than the direct product decomposition given by Wedderburn's theorem, the formula will not be entirely nice. 

Recall that in $\mathbb{C}$, we may write the standard branch of the logarithm as $Log(z) = log(|z|) + i \mathrm{Arg}(z)$, where in this case $|z|$ is just the standard complex modulus or norm $|z| = |x+iy| = \sqrt{x^2 + y^2}$, and $Arg$ is the function which takes a complex number and returns the standard angle of that number in polar form.

To summarize for the reader less familiar with complex analysis, given a number $z \in \mathbb{C}$ there exists $\rho \in \mathbb{R}^+$ and $\theta \in \mathbb{R}$ such that $z = \rho e^{i \theta}$, this being called the polar form of the complex number $z$. This is easy to see geometrically, given Euler's formula $e^{ix} = \cos(x) + i \sin(x)$, we see that $e^{i \theta}$ parameterizes the unit circle, and thus as we vary $\rho$ and $\theta$ we can reach any complex number. Furthermore, if $\rho \neq 0$, and $\theta \in [0,2\pi)$, $\rho$ and $\theta$ are unique, and hence we may define the standard argument $\mathrm{Arg}(z)$ by returning this unique $\theta$ given a number $z = \rho e^{i \theta}$.

Let $\mathcal{A}$ be a principal algebra -- that is, an algebra which may presented as $\mathbb{R}[k]/\langle p(k) \rangle$ for some $p(x) \in \mathbb{R}[x]$, and consider $z \in \mathrm{Ld}(\mathcal{A})$, then for some suggestively named constants $\ln(\rho) \in \mathbb{R}^+$ and  $\theta_1, \dots , \theta_{n-1} \in \mathbb{R}$, we have $z = e^{\ln(\rho) + k \theta_1 + \dots + k^{n-1} \theta_{n-1}} = e^{\ln(\rho)}e^{k \theta_1 + \dots + k^{n-1} \theta_{n-1}} = \rho e^{k \theta_1 + \dots + k^{n-1} \theta_{n-1}} $. This motivates the following definition:

\begin{definition}
Given a principal algebra $\mathcal{A}$ with generator $k$ and an element $z \in \mathrm{Ld}(\mathcal{A})$, then $z = e^{x_1 + k \theta_1 + \dots + k^{n-1} \theta_{n-1}}$ for some $x_1, \theta_1, \dots, \theta_{n-1} \in \mathbb{R}$. We say that $z$ is in generalized polar form when it is written as $ z = \rho e^{k \theta_1 + \dots + k^{n-1} \theta_{n-1}}$, where $\rho = e^{x_1}$ is an arbitrary positive real number.
\end{definition}

Notice however, that there is nothing in the above definition regarding uniqueness of the constants $\rho, \theta_1, \dots , \theta_{n-1}$. To show the uniqueness of $\rho$, and to provide a specific formula, we must extend the results of Freese \cite{Freese15}.

\begin{definition}
Given an algebra $\mathcal{A}$ with basis $\beta$, we define the Pythagorean function on the algebra $F(z)$ by $F(z) = \det(M_\beta(z))$. 
\end{definition}

\noindent Given this, Freese \cite{Freese15} proves the so-called k-thagorean theorem for the algebras $\mathcal{H}_k$ and $\mathcal{C}_k$:

\begin{theorem}
Let $j$ be the generator of the algebra $\mathcal{H}_k$, then for all $\theta \in \mathbb{R}$, $F(e^{j \theta}) = 1$. Similarly, if $i$ is the generator for the algebra $\mathbb{C}_k$, then $F(e^{i \theta}) = 1$
\end{theorem}

For example, this gives us in the familiar two dimensional cases the identity $\sin^2(\theta) + \cos^2(\theta) = 1$ by applying the theorem to $\mathbb{C}$ and $\cosh^2(\theta) - \sinh^2(\theta) = 1$ by applying it to $\mathcal{H}$.

We would like to extend this result to arbitrary principal algebras $\mathcal{A} = \mathbb{R}[k]/\langle p(k) \rangle$ where $p(k) = k^n + a_{n-1} k^{n-1} + \dots + a_1 k + a_0$. For reasons to be made clear in the following lemma, we apply the isomorphism $k \mapsto k - \frac{a_{n-1}}{n}$ so that $\mathcal{A} \cong \mathbb{R}[k]/\langle p(k - \frac{a_{n-1}}{n}) \rangle $, where now the polynomial $p(k - \frac{a_{n-1}}{n})$ is depressed. In other words, $p(k - \frac{a_{n-1}}{n}) = k^n + c_{n-2}k^{n-2} + \dots + c_1 k + c_0$ for some constants $c_0, c_1, \dots, c_{n-2} \in \mathbb{R}$. Thus, without loss of generality, we suppose $p(k)$ is in fact a depressed polynomial.

Before we move on, it is also important to note that every semisimple algebra has a presentation of the form outlined in the above paragraph, since we may simply choose $n$ coprime first order irreducible polynomials $p_1, \dots, p_n$ and $m$ second order irreducible polynomials $q_1, \dots, q_m$ and let $p(k) = p_1(k)p_2(k)\dots p_n(k)q_1(k)q_2(k)\dots q_m(k)$ so that by the Chinese remainder theorem $\mathbb{R}[k]/\langle p(k) \rangle \cong \mathbb{R}^n \times \mathbb{C}^m$.\footnote{For a more in depth discussion, see \cite{algebrapaper}.}

\begin{lemma}
\label{lem:generalized_k_thagorean_lemma}
Let $\mathcal{A}$ be a semisimple principal algebra, then $\mathcal{A} \cong \mathbb{R}[k]/\langle p(k) \rangle$ for some monic depressed $p(k) \in \mathbb{R}[k]$, where $p(k) = k^n + c_{n-2} k^{n-2} + \dots + c_1 k + c_0$. Given such a presentation, and basis $\beta = {1, k, \dots, k^{n-1} }$ consider $e^{k \theta} : \mathbb{R} \rightarrow \mathbb{R}$, and suppose $M_\beta(e^{k \theta}) = [u_1(\theta) | u_2(\theta) | \dots | u_n(\theta)]$, then the derivative of each column $u_i(\theta)$ is a linear combination of the remaining columns.

\begin{proof}
By theorem \ref{thm:regular_rep_cols} we know that $u_2 = k u_1, u_3 = k^2 u_1, \dots, u_n = k^{n-1} u_1$, so $e^{k \theta} = [u_1 | k u_1 | \dots | k^{n-1} u_1]$. Also, we have both 

\begin{equation*}
\begin{aligned}
\frac{d}{d\theta} e^{k z} = k e^{k \theta} &= k [u_1 | k u_1 | \dots | k^{n-1} u_1] \\ &= [k u_1 | k^2 u_1 | \dots | k^n u_1] \\ &= [u_2 | u_3 | \dots | k^n u_1]
\end{aligned}
\end{equation*}
and 
\begin{equation*}
\frac{d}{d\theta} e^{k \theta} = [u_1' | u_2' | \dots | u_n']
\end{equation*}
so by equating columns we obtain: $$ u_1' = u_2, u_2' = u_3, \dots , u_{n-1}' = u_n $$ So it remains to show the claim for $u_n$. Recall that $k^n + c_{n-2} k^{n-2} + \dots + c_1 k + c_0 = 0$ in our algebra, so $k^n = - c_{n-2} k^{n-2} - \dots - c_1 k - c_0 $, and hence: 

\begin{equation*}
\begin{aligned}
c_n' = k^n u_1 &= (- c_{n-2} k^{n-2} - \dots - c_1 k - c_0) u_1 \\ &= -c_{n-2} k^{n-2} u_1 - \dots - c_1 k u_1 - c_0 u_1 \\ &= -c_{n-2} u_{n-1} - \dots - c_1 u_2 - c_0 u_1
\end{aligned}
\end{equation*}
which is a linear combination of the other columns, thus completing the proof.
\end{proof}
\end{lemma}

\noindent From this lemma, we may derive the direct generalization of Freese's k-thagorean theorem: 

\begin{theorem}
Let $\mathcal{A} = \mathbb{R}[k]/\langle p(k) \rangle$ where $p(k)$ is depressed and monic, then for any $\theta \in \mathbb{R}$, $F(e^{k \theta}) = 1$.

\begin{proof}
Let $M(\theta) = M_\beta(e^{k \theta})$, then by definition of the determinant and applying the generalized product rule termwise we have:
\begin{equation*}
\begin{aligned}
\frac{d F(e^{k \theta})}{d\theta} &= \frac{d}{d\theta} \sum \epsilon_{i_1 \dots i_n} M_{i_1 1} M_{i_2 2} \dots M_{i_n n} \\ 
&= \sum \epsilon_{i_1 \dots i_n} \frac{d M_{i_1 1}}{d\theta} M_{i_2 2} \dots M_{i_n n} + \dots + \epsilon_{i_1 \dots i_n} M_{i_1 1} M_{i_2 2} \dots \frac{d M_{i_n n}}{d\theta} \\
 &= \sum \epsilon_{i_1 \dots i_n} \frac{d M_{i_1 1}}{d\theta} M_{i_2 2} \dots M_{i_n n} + \dots + \sum \epsilon_{i_1 \dots i_n} M_{i_1 1} M_{i_2 2} \dots \frac{d M_{i_{n-1} n-1}}{d\theta} M_{i_n n} \\ &- \sum  \sum_{j=0}^{n-2} c_j \epsilon_{i_1 \dots i_n} M_{i_1 1} M_{i_2 2} \dots M_{i_j j} \\
 &= \det[M_2 | M_2 | \dots | M_n] + \det [M_1 | M_3 | M_3 | \dots | M_n]  + \dots \\ &+ \det [M_1 | \dots | M_{n-1} | M_{n-1} | M_n ]  - \sum_{j=0}^{n-2} c_j \det[M_1 | \dots | M_j] \\
 &= 0
\end{aligned}
\end{equation*}
since by the previous lemma, working now at the component level of the matrix instead of at the column level, we obtained $\frac{d M_{i ,k}}{d\theta} = M_{i+1, k}$ for $i = 1, \dots, n-1$ and $$\frac{d M_{n,k}}{d\theta} = -\sum_{j = 0}^{n-2} c_n M_{j+1,k}$$ for $i = n$. Hence, after applying the generalized product rule and splitting up the sum, as demonstrated explicitly in the calculation above, we obtain a sum of determinants with repeated columns, and hence, $$ \frac{d F(e^{k\theta})}{d\theta} = 0 $$ So $F(e^{k\theta})$ must be constant. Furthermore, $F(e^{k 0}) = F(1) = 1$, so for all $\theta \in \mathcal{A}$ we must have $F(e^{k\theta}) = 1$.
\end{proof}
\end{theorem}

Notice that this theorem reproduces the usual k-thagorean theorem when $\mathcal{A} = \mathcal{H}_n$ or $\mathcal{A} = \mathbb{C}_n$. It essentially says that in some sense elements $e^{k \theta}$ are ``units'' of the algebra, since $F(e^{k \theta}) = 1$, although $F$ is not in general a norm, we would like to extend this statement to give a ``size" to all elements $z \in \mathrm{Ld}(\mathcal{A})$ with respect to the Pythagorean function $F$ which reproduces the radius of an algebra element with respect to the generalized polar form. Notice that if $z = \rho e^{k \theta_1 + k^2 \theta_2 + \dots + k^{n-1} \theta_{n-1}}$ and $F(e^{k \theta_1 + k^2 \theta_2 + \dots + k^{n-1} \theta_{n-1}}) = 1$ for all $\theta_1, \theta_2, \dots, \theta_{n-1} \in \mathbb{R}$ then $F(\rho e^{k \theta_1 + k^2 \theta_2 + \dots + k^{n-1} \theta_{n-1}}) = \rho^n F(e^{k \theta_1 + k^2 \theta_2 + \dots + k^{n-1} \theta_{n-1}}) = \rho^n$. Hence, we obtain a formula for $\rho$ by setting $\rho = \sqrt[n]{F(z)}$. However, as we will see in the following theorem, it does not seem to be the case that $F(e^{k \theta_1 + k^2 \theta_2 + \dots + k^{n-1} \theta_{n-1}}) = 1$ for all $\theta_1, \theta_2, \dots, \theta_{n-1} \in \mathbb{R}$ in general:

\begin{proposition}
\label{prop:only_hn_cn}
Given an $n$-dimensional principal algebra with generator $k$ such that $k^n = a_{n-1} k^{n-1} + \dots + a_1 k + a_0$, the derivatives of each of the columns of the functions $F(e^{k \theta}), F(e^{k^2 \theta}), \dots, F(e^{k^{n-1} \theta})$ are multiples of the remaining columns if and only if $a_{n-1} = a_{n-2} = \dots = a_1 = 0$.
\begin{proof}
Recall from lemma \ref{lem:generalized_k_thagorean_lemma} that in order for the derivatives of each of the columns of $e^{k \theta}$ to be a linear combination of the other columns, we need $a_{n-1} = 0$, where $k^n = a_{n-1} k^{n-1} + \dots + a_1 k + a_0$. Following a similar argument, let $e^{k^2 \theta} = [u_1 | u_2 |\dots| u_n] = [u_1|k u_1 | \dots | k^{n-1} u_1]$, then $$\frac{d}{d\theta} e^{k^2 \theta} = k^2 e^{k^2 \theta} = [k^2 u_1 | k^3 u_1| \dots | k^{n} u_1|k^{n+1} u_1] $$ and hence $u_{n-1}' = k^n u_1 = a_{n-1} k^{n-1} u_1 + a_{n-2} k^{n-2} u_1 + \dots + a_1 k u_1 + a_0 u_1 = a_{n-1} u_n + a_{n-2} u_{n-1} + \dots + a_1 u_2 + a_0 u_1$. Therefore, in order for $u_{n-1}'$ to be a linear combination of the other columns we need $a_{n-2} = 0$.

Applying the same technique to $e^{k^3 \theta}, \dots, e^{k^{n-1} \theta}$, we find that $a_{n-2}, a_{n-3}, \dots, a_2, a_1$ must all be zero. Furthermore, simple calculations also show that this condition is also sufficient for the derivatives of each of the columns of the functions $F(e^{k \theta}), F(e^{k^2 \theta})$, $\dots$, $F(e^{k^{n-1} \theta_{n-1}})$ to be multiples of the remaining columns. Therefore, this is the case if and only if we are working over an algebra with $a_{n-1} = a_{n-2} = \dots = a_1 = 0$.
\end{proof}
\end{proposition}

\noindent Hence, although this is not a definitive result, it seems unlikely that $$F(e^{k \theta_1 + k^2 \theta_2 + \dots + k^{n-1} \theta_{n-1}}) = 1$$ for all $\theta_1, \theta_2, \dots, \theta_{n-1} \in \mathbb{R}$ for general principal algebras, though still technically possible. Given the preconditions of proposition \ref{prop:only_hn_cn}, this is true for $\Hcal_n, \Ccal_n,$ and $\mathbf{\Gamma}_n$. We leave it as an open question whether or not this is possible for other principal algebras, but again, given proposition \ref{prop:only_hn_cn}, this seems unlikely. Although, there may be other fruitful approaches to this problem, for example, considering bases other than principal bases (i.e. bases of the form $\{ 1, k, \dots, k^{n-1} \}$).

Thus, at least for the moment, we develop the theory of the modulus for an algebra restricted to the case where $\Acal = \Ccal_n, \Hcal_n$ or $\mathbf{\Gamma}_n$.

\begin{corollary}
\label{cor:modulus}
Let $\mathcal{A} = \Ccal_n, \Hcal_n$, or $\mathbf{\Gamma}_n$. In other words, $\Acal$ is a principal algebra with generator $k$ such that $k^n = a$ for some $a \in \mathbb{R}$. Given a number $z \in \mathrm{Ld}(\mathcal{A})$, let $z = \rho e^{k \theta_1 + \dots + k^{n-1} \theta_{n-1}}$ be a representation of $z$ in generalized polar form. Then $\rho$ is the unique positive real number that satisfies the equation $F(z) = \rho^n$, and hence is given by $\rho = \sqrt[n]{F(z)}$. Thus, given $\mathcal{A}$ has basis $\beta = \{1,k,\dots,k^{n-1}\}$, we define the $\mathcal{A}$-modulus $|z|_\beta = \sqrt[n]{F(z)}$, which we will often write as simply $|z|$ if the basis is understood.

\begin{proof}
Given an algebra such that $k^n = a$, from proposition \ref{prop:only_hn_cn} it follows that $F(e^{k \theta_1}) = F(e^{k^2 \theta_2}) = \dots = F(e^{k^{n-1} \theta_{n-1}}) = 1$ and hence, since $F = \det \comp M$, and both $\det$ and $M$ satisfy the homomorphism property, $F$ also satisfies the homomorphism property, we obtain: 
\begin{equation*}
\begin{aligned}
&\;\;\;\;\; F(e^{k \theta_1}) F(e^{k^2 \theta_2}) \dots F(e^{k^{n-1} \theta_{n-1}}) \\ &= F(e^{k \theta_1}e^{k^2 \theta_2} \dots e^{k^{n-1} \theta_{n-1}}) \\ &= F(e^{k \theta_1 + k^2 \theta_2 + \dots + k^{n-1} \theta_{n-1}}) \\ &= 1
\end{aligned}
\end{equation*}

\noindent Now, let $z = \rho e^{k \theta_1 + k^2 \theta_2 + \dots + k^{n-1} \theta_{n-1}}$, then: 
\begin{equation*}
\begin{aligned}
F(z) &= F(\rho e^{k \theta_1 + k^2 \theta_2 + \dots + k^{n-1} \theta_{n-1}}) \\ &= \det(M_\beta(\rho e^{k \theta_1 + k^2 \theta_2 + \dots + k^{n-1} \theta_{n-1}})) \\ &= \det(\rho M_\beta(e^{k \theta_1 + k^2 \theta_2 + \dots + k^{n-1} \theta_{n-1}})) \\ &= \rho^n \det(M_\beta(e^{k \theta_1 + k^2 \theta_2 + \dots + k^{n-1} \theta_{n-1}})) \\ &= \rho^n
\end{aligned}
\end{equation*}

\noindent Hence, since $\rho \geq 0$, we obtain that $\rho = \sqrt[n]{F(z)}$ is unique.

\end{proof}
\end{corollary}

\begin{theorem}
Let $Log(z)$ be any inverse to the exponential function in an algebra $\Acal = \Ccal_n, \Hcal_n,$ or $\mathbf{\Gamma}_n$, then $\mathrm{Re}(Log(z)) = \log(|z|)$.

\begin{proof}
let $z = e^{x_1 + k x_2 + \dots + k^{n-1} x_n} \in \mathrm{Ld}(\mathcal{A})$, then recalling from the proof of corollary \ref{cor:modulus} that $F(e^{k x_2 + \dots + k^{n-1} x_n}) = 1$, and hence that $|e^{k x_2 + \dots + k^{n-1} x_n}| = 1$, we obtain: 

\begin{equation*}
\begin{aligned}
\log(|z|) &= \log(|e^{x_1 + k x_2 + \dots + k^{n-1} x_n}|) \\ &= \log(|e^{x_1} e^{k x_2 + \dots + k^{n-1} x_n}|) \\ &= \log(e^{x_1} |e^{k x_2 + \dots + k^{n-1} x_n}|) \\ &= \log(e^{x_1}) \\ &= x_1
\end{aligned}
\end{equation*}

\noindent And hence, $\log(|z|)$ is the real part of the function $\mathrm{Log}(z)$
\end{proof}
\end{theorem}

Thus, we have found one half of the analogue of the formula from complex analysis that $\mathrm{Log}(z) = \log(|z|) + i \mathrm{Arg}(z)$. Unfortunately, in more general algebras the notion of argument is somewhat complicated, at least from a formulaic perspective.

Restricting our attention to the two dimensional case, let $z \in \mathbb{C}$ with polar form $z = \rho e^{i \theta} = \rho(\cos(\theta) + i \sin(\theta))$. To find the function $\mathrm{Arg}$ returning the unique $\theta \in [0,2\pi)$ such that $z = \rho e^{i \theta}$, so long as $\cos(\theta) \neq 0$, if $x = \mathrm{Re}(z)$ and $y = \mathrm{Im}(z)$ we have $\frac{y}{x} = tan(\theta) \implies \theta = \mathrm{Arg}(z) = \tan^{-1}(\frac{y}{x})$, giving us a formula for the argument.

Similarly, for $z = \rho e^{j \psi} = \rho (\cosh(\psi) + j \sinh(\psi))$, we get an even better formula for the argument, since we have $\cosh(x) \neq 0$ for all $x \in \mathbb{R}$. Hence, since $x = \mathrm{Re}(z) = \rho \cosh(\psi)$, and $y = \mathrm{Im}(z) = \rho \sinh(\psi)$, we obtain $\tanh(\psi) = \frac{y}{x} \implies \psi = \mathrm{Arg}(z) = \tanh^{-1}(\frac{y}{x})$ as a formula for the argument of all of $\mathrm{Ld}(\mathcal{H})$.

It is a natural question then to ask, and part of the motivation for this section of the paper, whether or not we can find similar formulas for the argument (i.e. the function giving the imaginary part of the logarithm for an algebra) in general. In particular, we hoped that in the case of $\mathcal{H}_n$ and $\mathbb{C}_n$, we could use the component functions of $e^{j \theta}$ and $e^{i \theta}$ -- what \cite{Freese15} calls the $k$-hyperbolic and $k$-trigonometric functions respectively, would allow us to find similar formulas for the argument on at least part of $\mathrm{Ld}(\mathcal{A})$.  Unfortunately, as soon as we move away from two dimensional algebras, the situation is complicated by the fact that the polar form of a number $z \in \mathrm{Ld}(\mathcal{A})$ is no longer simply of the form $\rho e^{k \theta}$, but is in general $\rho e^{k \theta_1 + \dots + k^{n-1} \theta_{n-1}} = \rho e^{k \theta_1} e^{k^2 \theta_2} \dots e^{k^{n-1} \theta_{n-1}}$, and thus the components of the polar form representation are now sums of products of the component functions of $e^{k \theta_1}, \dots, e^{k^{n-1} \theta_{n-1}}$.

\begin{example}
In $\mathcal{H}_3$, making use of the identities $\cosh_3(jz) = \cosh_3(z), \sinh_{3,1}(jz) = j \sinh_{3,1}(z)$, and $\sinh_{3,2}(jz) = j^2 \sinh_{3,2}(z)$ we have: 
\begin{equation*}
\begin{aligned}
& e^{j \theta + j^2 \psi} = e^{j \theta} e^{j^2 \psi} \\ &= (\cosh_3(\theta) + j \sinh_{3,1}(\theta) + j^2 \sinh_{3,2}(\theta))(\cosh_3(j \psi) + j \sinh_{3,1}(j \psi) + j^2 \sinh_{3,2}(j \psi )) \\ &= (\cosh_3(\theta) + j \sinh_{3,1}(\theta) + j^2 \sinh_{3,2}(\theta))(\cosh_3(\psi) + j^2 \sinh_{3,1}(\psi) + j \sinh_{3,2}(\psi )) \\ &= (\cosh_3(\theta)\cosh_3(\psi) + \sinh_{3,1}(\theta)\sinh_{3,1}(\psi) + \sinh_{3,2}(\theta)\sinh_{3,2}(\psi)) \\ &+ j(\cosh_3(\theta)\sinh_{3,2}(\psi) + \sinh_{3,1}(\theta)\cosh_3(\psi) + \sinh_{3,2}(\theta)\sinh_{3,1}(\psi)) \\ &+ j^2(\cosh_3(\theta)\sinh_{3,1}(\psi) + \sinh_{3,2}(\theta)\cosh_3(\psi) + \sinh_{3,1}(\theta)\sinh_{3,2}(\psi) )
\end{aligned}
\end{equation*}

\end{example}

And thus, we leave the question of if, or to what extent, nice formulas for the standard argument function may be found for higher dimensional algebras as an open question, and for now define $\mathrm{Arg}(z)$ simply by passing the logarithm on $\mathbb{R}^m \times \mathbb{C}^k$ through the isomorphism to a principle algebra $\mathcal{A}$.

\begin{definition}
Let $\mathcal{A}$ be a principle semisimple algebra, and let $Log(z)$ be the logarithm induced by the isomorphism $\phi$ from $\mathcal{A}$ to $\mathbb{R}^n \times \mathbb{C}^m$ given a particular choice of branch cuts $B_1, \dots, B_m$ for each of the complex components, then given $w \in \Ld(\Acal)$, the branch cut of the logarithm selects the unique element $z = x_1 + k \theta_1 + \dots + k^{n-1} \theta_{n-1} \in \phi^{-1}(\mathbb{R}^n \times B_1 \times \dots \times B_m)$ such that $e^z = w$. Thus, we define $\mathrm{Arg}(z) = k \theta_1 + k^2 \theta_2 + \dots + k^{n-1} \theta_{n-1}$ for all $z \in \mathrm{Ld}(\mathcal{A})$.
\end{definition}

And hence, since we have already shown explicitly that the first component of $Log_\mathcal{A}(z)$ is always $\log(|z|)$ if $\Acal = \Hcal_n, \Ccal_n$ or $\mathbf{\Gamma}_n$ it follows that:

\begin{corollary}
$Log_\mathcal{A}(z) = log(|z|) + Arg(z)$ if $\Acal = \Hcal_n, \Ccal_n$, or $\mathbf{\Gamma}_n$.
\end{corollary}

\subsection{Generalized Trigonometry and Open Problems}

We conclude our section on generalized polar form by generalizing some of the results of classical trigonometry to the component functions of generalized complex exponentials $e^{k z}$ for $k$ the generator of some semisimple principal algebra, in addition to working out specific formulas for some low dimensional algebras. Besides being interesting in their own right, we believe that these results may be useful in future researchers who wish to tackle the problem of finding explicit formulas for $\mathrm{Arg}(z)$ in terms of the component functions of $e^{k z}$.

\begin{definition}
Given a principle semisimple algebra $\mathcal{A} = \mathbb{R}[k]/\langle p(k) \rangle$, we say the component functions $s_1, \dots, s_n$ of the generalized complex exponential function $e^{kz} = s_1 + k s_2 + \dots + k^{n-1} s_n$ are the generalized trigonometric functions of the algebra $\mathcal{A}$.

In particular, in the case that $\mathcal{A} = \mathcal{H}_n$, we use the naming scheme $s_1 = \cosh_n(z)$, and $s_2 = \sinh_{n,1}(z), s_3 = \sinh_{n,2}(z), \dots, s_n = \sinh_{n,n-1}(z)$. Similarly, in the case that $\mathcal{A} = \mathbb{C}_n$, we use the naming scheme $s_1 = \cos_n(z)$, and $s_2 = \sin_{n,1}(z), s_3 = \sin_{n,2}(z), \dots, s_n = \sin_{n,n-1}(z)$.
\end{definition}

Thus, our generalized trigonometric functions generalize what Freese calls the $k$-trigonometric and the $k$-hyperbolic functions; themselves generalizations of $\sin,\cos,\sinh$ and $\cosh$.\footnote{\label{footnote:A-ODEs}In fact, these generalized trigonometric functions are again generalized by a class of canonical solutions to particular ODEs over an algebra, which we have studied in am unpublished paper on the subject.} Our goal now is to study the properties of these functions, and their relations to each other, generalizing commonly known results about $\sin, \cos, \sinh$ and $\cosh$.

\begin{proposition}[Generalized Adding Angles Formula]
\label{prop:generalized_adding_angle}
Let $\mathcal{A}$ be a principal algebra with generator $k$, then consider $e^{kz} = s_1(z) + k s_2(z) + \dots + k^{n-1} s_n(z)$. Using laws of exponents we obtain $e^{k(\alpha + \beta)} = e^{k \alpha} e^{k \beta}$. Hence: $$ \sum_{i = 1}^{n} s_i(\alpha + \beta) k^{i-1} = (\sum_{j=1}^n s_j(\alpha) k^{j-1})(\sum_{l=1}^n s_l(\beta) k^{l-1}) = \sum_{j,l = 1}^{n} s_j(\alpha) s_l(\beta) k^{j+l-2} $$ and thus by equating coefficients of $1, k, \dots, k^{n-1}$ we obtain formulas for $s_i(\alpha + \beta)$ in terms of products, sums, and scalar multiples of functions of the form $s_j(\alpha)$ and $s_l(\beta)$. We call these the generalized adding angle formulas for the algebra $\mathcal{A}$.
\end{proposition}

\begin{example}
\label{ex:gen_add_angle}

Following the procedure laid out in proposition \ref{prop:generalized_adding_angle} to find the generalized adding angles formula for $\mathcal{H}_3$ we find: $$ e^{j(\alpha + \beta)} = \cosh_3(\alpha + \beta) + j \sinh_{3,1}(\alpha + \beta) + j^2 \sinh_{3,2}(\alpha + \beta) $$ and

\begin{equation*}
\begin{aligned}
e^{j\alpha}e^{j\beta} &= (\cosh_3(\alpha) + j \sinh_{3,1}(\alpha)+j^2\sinh_{3,2}(\alpha))(\cosh_3(\beta) + j \sinh_{3,1}(\beta)+j^2\sinh_{3,2}(\beta)) \\ & = (\cosh_3(\alpha)\cosh_3(\beta) + \sinh_{31}(\alpha)\sinh_{3,2}(\beta) + \sinh_{3,2}(\alpha) \sinh_{3,1}(\beta)) \\ &+ j(\sinh_{3,1}(\alpha)\cosh_3(\beta) + \cosh_3(\alpha) \sinh_{3,1}(\beta) + \sinh_{3,2}(\alpha) \sinh_{3,2}(\beta)) \\ &+ j^2(\sinh_{3,2}(\alpha) \cosh_3(\beta) + \cosh_3(\alpha) \sinh_{3,2}(\beta) + \sinh_{3,1}(\alpha) \sinh_{3,1}(\beta))
\end{aligned}
\end{equation*}

\noindent And hence, by equating coefficients of $e^{j(\alpha + \beta)}$ and $e^{j \alpha} e^{j \beta}$ we obtain the formulas: $$\cosh_3(\alpha + \beta) = \cosh_3(\alpha) \cosh_3(\beta) + \sinh_{3,1}(\alpha) \sinh_{3,2}(\beta) + \sinh_{3,2}(\alpha) \sinh_{3,1}(\beta)$$ $$\sinh_{3,1}(\alpha + \beta) = \sinh_{3,1}(\alpha)\cosh_3(\beta) + \cosh_3(\alpha) \sinh_{3,1}(\beta) + \sinh_{3,2}(\alpha) \sinh_{3,2}(\beta)$$ $$ \sinh_{3,2}(\alpha + \beta) = \sinh_{3,2}(\alpha) \cosh_3(\beta) + \cosh_3(\alpha) \sinh_{3,2}(\beta) + \sinh_{3,1}(\alpha) \sinh_{3,1}(\beta)$$ 

\noindent For all $\alpha, \beta \in \mathcal{H}$.

\end{example}

\begin{proposition}[Generalized De Moivre's formulas]
\label{prop:generalized_demoivre}
Let $\mathcal{A}$ be a principal algebra with generator $k$, then for $l \in \mathbb{N}$, consider the equation $e^{l\alpha} = {e^\alpha}^l$. At the level of component functions then we have
\begin{equation*}
\begin{aligned}
& (\sum_{i_1=1}^l s_{i_1} k^{i_1-1})(\sum_{i_2=1}^l s_{i_2} k^{i_2-1})\dots(\sum_{i_n=1}^l s_{n_1} k^{i_n-1}) \\ &= \sum_{i_1,\dots,i_n = 1}^l s_{i_1}(\alpha)\dots s_{i_n}(\alpha) k^{i_1 + \dots + i_n - n} \\ &= \sum_{j=1}^{l} s_j(n\alpha) k^{j-1}
\end{aligned}
\end{equation*}
Hence, by equating the coefficients for $1,k,\dots,k^{n-1}$ we obtain equations relating $s_i(n\alpha)$ to sums of products and scalar multiples of $s_1(\alpha), s_2(\alpha), \dots , s_n(\alpha)$. We call these the generalized De Moivre's formulas.
\end{proposition}

\begin{example}
In $\mathcal{H}$, following the procedure in proposition \ref{prop:generalized_demoivre} we obtain:

\begin{equation*}
\begin{aligned}
e^{j \psi} &= \cosh(3\psi) + j \sinh(3 \psi) \\ &= (\cosh(\psi) + j \sinh(\psi))^3 \\ &= \cosh^3(\psi) + \cosh(\psi)\sinh^2(\psi) + 2\cosh(\psi)\sinh^2(\psi) \\ &+ j(\sinh^3(\psi) + \sinh(\psi)\cosh^2(\psi) + 2 \sinh(\psi)\cosh^2(\psi)) 
\end{aligned}
\end{equation*}

\noindent And hence, equating coefficients we obtain: $$ \cosh(3\psi) = \cosh^3(\psi) + \cosh(\psi)\sinh^2(\psi) + 2\cosh(\psi)\sinh^2(\psi) $$ $$ \sinh(3 \psi) = \sinh^3(\psi) + \sinh(\psi)\cosh^2(\psi) + 2 \sinh(\psi)\cosh^2(\psi) $$
\end{example}

Beyond these formulas, we suspect that there are other interesting questions that can be asked about these functions. For example, in $\mathcal{H}_n$ it is simple to see from the series expansions of $\cosh_n(z), \sinh_{n,1}(z), \dots, \sinh_{n,n-1}(z)$ that $\cosh_n(jz) = \cosh_n(z)$ and that $\sinh_{n,k}(jz) = j^k \sinh(z)$, with similar results that can be derived for $\mathcal{C}_n$. Thus, we wonder:

\begin{open}
Given an arbitrary principal algebra $\mathcal{A}$ with generator $k$ and generalized trigonometric functions $s_1, \dots, s_n$, are there any interesting relations between $s_i(kz)$ and the other generalized trigonometric functions?
\end{open}

While we don't expect there to be very nice relations between $s_i(kz)$ and $s_1, \dots, s_n$ in general, after all, as we saw in example \ref{ex:gen_add_angle}, even for $\mathcal{H}_3$ these generalized trigonometric formulas can be complex, a more interesting question for further research into the properties of these generalized trigonometric functions might be:

\begin{open}
Are $\sin$ and $\cos$ the only periodic generalized trigonometric functions?
\end{open}

Starting only from the series definition of $\sin$ and $\cos$, proving the periodicity of the two functions is somewhat involved. Notably however, the proof involves the Pythagorean theorem. Thus, we wonder how the generalized $k$-thagorean theorem might play a role in the proof of periodicity or lack thereof for the generalized trigonometric functions $s_1, \dots, s_n$. 

Such a general result on the periodicity of generalized trigonometric functions may even provide interesting insights into the solution sets of linear ODEs, as we noted in footnote \ref{footnote:A-ODEs}, generalized trigonometric functions may be seen as canonical solutions to certain ODEs over an algebra.

\pagebreak
\bibliography{sources}

\end{document}